\magnification\magstep1

 \baselineskip = 18pt
 \def\n{\noindent}

\def\adh{\overline}
\def \rat{ {\rm Q}\kern-.65em {}^{{}_/ }}

 \magnification\magstep1
 \baselineskip = 18pt
 \def\n{\noindent}
 \overfullrule = 0pt
 \def\qed{{\hfill{\vrule height7pt width7pt
depth0pt}\par\bigskip}}
 \centerline{\bf On the ``local theory'' of operator spaces}\medskip
 
  \centerline{{by Gilles Pisier}\footnote*{Supported in
part by N.S.F. grant DMS 9003550}}\bigskip
 
 In Banach space theory, the ``local theory'' refers to the collection of
 finite dimensional methods and ideas which are used to study infinite
 dimensional spaces (see e.g. [P4,TJ]). It is natural to
try to develop an
 analogous theory in the recently developed category of operator spaces
 [BP,B1-2,BS,ER1-7,Ru]. The object of this paper is to
start such a theory. We plan to present a more thorough
discussion of the associated tensor norms in a future
publication.
 
We refer to [BP,B1-2, ER1-7] for the definition and
the main properties of operator spaces. We merely recall
that an operator space is a Banach space isometrically
embedded into the space $B(H)$ of all bounded operators
on a  Hilbert space $H$, and that in the category of
operator spaces, the morphisms are the completely bounded
maps (in short cb) for which we refer the reader to [Pa1].
If $E,F$ are operator spaces, we denote by
$E\otimes_{\min}F$ their minimal (or spatial) tensor
product. We denote by $\adh E$ (or $\adh H$) the space 
$E$ (or $H$) equipped with the conjugate complex
multiplication.  Note that $\overline E^*$ can be identified with the
antidual
 of $E$ and the elements of $(E\otimes \overline E)^*$ can be viewed as
 sesquilinear forms on $E\times E$.

  Recently [P1-3] we introduced the
analogue of Hilbert space in the category of operator
spaces. We proved that there is a Hilbert space ${\cal H}$
and a sequence of operators $T_n\in B({\cal H})$ such that
for all finitely supported sequence $(a_n)$ in $B(\ell_2)$
we have  
$$\left\|\sum\nolimits  T_n \otimes
a_n\right\|_{B({\cal H}\otimes
 \ell_2)} = \left\|\sum  a_n\otimes \bar
a_n\right\|^{1/2}_{B(\ell_2
 \otimes\overline\ell_2)}.\leqno(1)$$
We denoted by $OH$ the closed span of $(T_n)$ and by
$OH_n$ the span of $T_1,...T_n$. We call $OH$ the
operator Hilbert space. For any operator $u\colon OH_n\to
E$ we have (cf.[P1-2])
$$\|u\|_{cb} =\|\sum_1^n u(T_i) \otimes
T_i\|_{E\otimes_{\min} OH_n}.$$

 Our main tool will be a variation (one more!) on the notion of 2-summing
 operator. Let $E$ be an operator space and let $Y$ be a Banach space.
 Following our previous work \break [P3] we will say that
an operator
 $u\colon
 \ E\to Y$ is $(2,oh)$-summing if there is a 
constant $C$ such that for all
 finite sequences $(x_i)$ in $E$ we have
 $$\sum \|u(x_i)\|^2 \le C^2\left\|\sum x_i\otimes \bar x_i\right\|
 _{E\otimes_{\rm min}\overline E}.\leqno(2)$$
 We will denote by $\pi_{2,oh}(u)$ the smallest constant $C$ for which this
 holds. Moreover for any integer $n$ we denote by
$\pi^n_{2,oh}(u)$ the smallest constant $C$ such that (2)
holds for all $n$-tuples $x_1,...,x_n$ in $E$. Recall
that the usual $2$-summing norm $\pi_{2}(u)$ of an
operator $u\colon E\to F$ between Banach spaces (resp. the
$2$-summing norm on $n$ vectors $\pi^n_{2,oh}(u)$)  is the
smallest constant $C$ such that 
  for all
 finite sequences  (resp. all $n$-tuples) $(x_i)$ in $E$ we
have
 $$\sum \|u(x_i)\|^2 \le C^2 \sup\{\sum|\xi(
x_i)|^2 |\ \xi\in E^* \ \|\xi\|\le 1\}.$$
Equivalently this means that (2) holds when $E$ is
embedded isometrically into a commutative $C^*$-subalgebra
of $B(H)$.  An
alternate definition of $\pi_{2}(u)$ (resp.  
$\pi^n_{2}(u)$)  is the
smallest constant $C$ such that 
 $$   \pi_2(uv) \leq C \leqno(3)$$
for all finite rank operators $v\colon \ell_2 \to E$
(resp. $v\colon \ell^n_2 \to E$) with $\|v\|\le 1.$ 
As observed in [P1], it
is easy to see using (1) that for every  bounded operator
$v\colon OH_n\to OH_n$ we have $\|v\|=\|v\|_{cb}$. It
follows that for any operator $u\colon OH\to E$
we have 
$$\pi_{2,oh}(u)=\pi_{2 }(u) \quad {\rm and }\quad
\pi^n_{2,oh}(u)=\pi^n_{2 }(u).$$
Similarly when $E$  is an operator space for any  $u\colon
E\to F$   the norm $\pi_{2,oh}(u)$
(resp. $\pi^n_{2,oh}(u)$) is the smallest constant $C$
such that  (3) holds for all finite rank $v\colon OH \to E$
(resp. $v\colon OH_n \to E$) with $\|v\|_{cb}\le 1.$ 
 
\n Since the cb-norm dominates the usual norm of
an operator $v\colon OH\to E$, it is easy to check
that for if $E$ is an operator space and $F$ a Banach
space
then every $2$-summing $u\colon E\to F$ is necessarily
$(2,oh)$-summing and we have
$$\pi_{2,oh}(u)\leq\pi_{2}(u).\leqno(4)$$

 \n By an important inequality
due to Tomczak-Jaegermann (see [TJ] p.143) we have for any 
rank $n$ operator  $u\colon
E\to F$ between Banach spaces $\pi_{2 }(u)\le
\sqrt 2 \pi^n_{2 }(u)$. This fact and the preceding
equalities yield that for any rank $n$ operator
$u\colon
E\to F$ between operator   spaces we have
$$\pi_{2,oh }(u)\le
\sqrt 2\pi^n_{2,oh }(u).\leqno (5)$$
  In [P3] the following
result (which is crucial for the present note) is
mentioned. Any operator $u\colon
 \ E\to OH$ (with domain an arbitrary operator space but
with range $OH$) which is $(2,oh)$-summing is necessarily
completely bounded and we have
$$\forall u\colon
 \ E\to OH\quad \|u\|_{cb} \leq \pi_{2,oh}(u).\leqno(6)$$

 \n An element $u$ in $E\otimes\overline E$ is called
positive if $u$ can be
 written as
 $$u = \sum^n_1 x_i\otimes \bar x_i \quad \hbox{with}\quad x_i\in E.$$
 In that case we will write $u\ge 0$. Equivalently this means that $\langle
 u, \xi\otimes \bar \xi\rangle \ge 0$ for all $\xi$ in $E^*$, so that $u\ge
 0$ iff $\langle u,v\rangle \ge 0$ for all $v\ge 0$ in $E^*\otimes
 \overline{E^*}$. More generally, a linear form $\varphi \in (E\otimes
 \overline E)^*$ will be called positive if $\varphi(x\otimes \bar x)\ge 0$
 for all $x$ in $E$. Note that this implies that $\varphi$ is symmetric,
 i.e. $\varphi$ satisfies $\varphi(x\otimes \bar y) = \varphi(y\otimes \bar
 x)$ (or equivalently $\varphi(u) \in {\bf R}$ for all {\it symmetric\/} $u$
 in $E\otimes \overline E$). We will denote by $K(E)$ the set of all the
 positive linear forms $\varphi$ in $(E\otimes\overline E)^*$ such that
 $$\sup\{\varphi(u) \mid u \in E\otimes \overline E,
u\ge 0, \|u\|_{E\otimes_{\rm min}\overline E}
 \le 1\} \le 1.$$
 Then it is rather easy to check that for all $u\geq 0$ in
$E\otimes \overline E$
 we have
 $$\|u\|_{E\otimes_{\rm min}\overline E} = \sup_{\varphi \in
K(E)}\varphi(u).\leqno (7)$$
 Indeed, consider $u = \sum\limits^n_1 x_i\otimes\bar y_i$
in $E\otimes
 \overline E$. Assume $E\subset B(H)$. Let $C_2$ be the space of all
 Hilbert-Schmidt operators on $H$, equipped with the Hilbert-Schmidt norm
 $\|\ \|_2$. Observe that for any $y$ in $C_2$ there is a
decomposition $y =
 a^+-a^- + i(b^+-b^-)$ with $a^+, a^-, b^+,b^-$ hermitian positive and such
 that
 $$\|a^+\|^2_2 + \|a^-\|^2_2 + \|b^+\|^2_2 + \|b^-\|^2_2 =
\|y\|^2_2.$$
 By definition of $E\otimes_{\rm min} \overline E$ we have
 $$\|u\| = \sup\left\{\left|{tr}\left(\sum x_iyy^*_iz\right)\right| \mid
 y,z\in C_2\  \|y\|_2\le 1,\  \|z\|_2 \le 1\right\}.$$
 Now assume $u\ge 0$, say $u=\sum\limits^n_1 x_i\otimes\bar
x_i$. Let $F(y,z)={tr}\left(\sum x_iyx^*_iz\right)$. Note
that $F(y,z)$ is positive when $y,z$ are both positive.
Then, by the decomposition recalled above the  
supremum of $F(y,z)$ when $y,z$ run over the unit ball
of $C_2$ is unchanged if
 we restrict it to positive operators $y,z$ in the unit ball of $C_2$. But
 if $y,z$ are positive in the unit ball of $C_2$ then 
the form defined by $\forall u = \sum\limits^n_1
x_i\otimes\bar y_i\in E\otimes \bar E$
 $$\eqalign{\varphi(u) &= {tr}\left(\sum x_iyy^*_iz\right) = {tr}\left(
 z^{1/2} x_iyy^*_iz^{1/2}\right)\cr
 &= \sum {tr}[(z^{1/2}x_iy^{1/2})
(z^{1/2}y_iy^{1/2})^*]}$$
 is clearly positive so that (7) follows.
 
 \proclaim Proposition 1. Let $E$ be an operator space, let $F\subset E$ be
 a closed subspace, and let $Y$ be a Banach space. Let $u\colon \ F\to Y$ be
 an operator and let $C$ be a constant. The following are
 equivalent.\medskip
 \item{(i)} $u$ is $(2,oh)$-summing with $\pi_{2,oh}(u) \le C$.
 \item{(ii)} There is a $\varphi$ in $K(E)$ such that
 $$\forall\ x \in F\qquad \|u(x)\|^2 \le C^2 \varphi(x\otimes \bar x).$$
 \item{(iii)} There is an extension $\tilde u\colon \ E\to Y$ such that
 $\tilde u_{|F} = u$ and
 $$\pi_{2,oh}(\tilde u) \le C.$$
 
 \noindent {\bf Proof:}\ Assume (i). Note that for all
$w\in F\otimes \overline F$ we have
$$\|w\|_{F\otimes_{\rm min}\overline
F}=\|w\|_{E\otimes_{\rm min}\overline E}.$$
 Hence by (7) we have
 $$\sum \|u(x_i)\|^2 \le C^2 \sup_{\varphi \in K(E)} \sum \varphi
 (x_i\otimes \bar x_i),$$
 for all finite sequences $x_i$ in $F$. By a classical application of the
 Hahn-Banach theorem it follows that there is a $\varphi$ in $K(E)$ such
 that
 $$\forall\ x\in E\qquad \|u(x)\|^2 \le C^2\varphi(x\otimes \bar x).$$
 (Indeed, one can reproduce the argument included e.g. in
[P4] p.~11 for
 2-summing operators and observe that $K(E)$ is convex so that the
 barycenter of a probability measure on $K(E)$ belongs to $K(E)$.) This
 proves (i)~$\Rightarrow$~(ii). Now assume (ii). Consider the scalar product
 $\langle x,y\rangle = \varphi(x\otimes \bar y)$ on $E$. Let us denote by
 $L^2(\varphi)$ the resulting Hilbert space (after
passing to the usual quotient and completing) and let
$J\colon \ E\to
 L^2(\varphi)$ be the natural inclusion. Observe that we trivially have by
 (7) $\pi_{2,oh}(J) \le 1$. We now introduce an operator
$v\colon \ J(F)\to
 Y$. For any element $y$ in $J(F)$ we can define if $y = J(x)$ with $x\in F$
 $$v(y) = u(x).$$
 Note that (ii) ensures that this definition is unambiguous and $\|v\| \le
 C$. Hence $v$ extends to an operator $v\colon \ \overline{J(F)} \to Y$ such
 that $\|v\| \le C$. Finally let $P$ be the orthogonal projection from
 $L^2(\varphi)$ onto $\overline{J(F)}$ and let $\tilde u = vPJ$. Clearly
 $\pi_{2,oh}(\tilde u) \le \|v\|\pi_{2,oh}(J) \le C$ and $\tilde u$ extends
 $u$. This proves (ii)~$\Rightarrow$~(iii). Finally (iii)~$\Rightarrow$~(i)
 is trivial. \qed
  
A fundamental inequality  in Banach space theory (originally due to Garling and
Gordon, see [P4,p.15]) says that
for any $n$-dimensional Banach space the identity operator
$I_E$ satisfies $\pi_{2}(I_E) = n^{1/2}.$ By (4) it
follows
that for any $n$-dimensional operator space
we have
$$\pi_{2,oh}(I_E) \leq n^{1/2}.$$
In that case the equality no longer holds, as shown by the
examples below.
However the following consequence of the upper bound still
holds in the category of operator spaces.
 \proclaim Theorem 2. Let $E$ be any $n$-dimensional operator space then
   there is an isomorphism
 $$u\colon \ E\to OH_n \quad \hbox{such that}\quad \pi_{2,oh}(u) = n^{1/2}$$
 and $\|u^{-1}\|_{cb} = 1$.

 \proclaim Corollary 3. For any $n$-dimensional operator
space  $E$ there are $n$  elements $x_1,...,x_n$ in $E$
such that 
$$\left\|\sum_1^n x_i\otimes \bar x_i\right\|
 _{E\otimes_{\rm min}\overline E} \le 1 \quad {\rm and
}\quad \sum_1^n \|x_i\|^2\geq \pi_{2,oh}(I_E)^2/2.$$

 \proclaim Corollary 4. For any $n$-dimensional $E$ there
is an isomorphism
 $u\colon \ OH_n\to E$ such that $\|u\|_{cb}\ \|u^{-1}\|_{cb} \le \sqrt n$.
 
 \proclaim Corollary 5. For any $n$-dimensional subspace
$E\subset B(H)$
 there is a projection $P\colon \ B(H) \to E$ such that
 $$ \|P\|_{cb} \le \sqrt n.$$

 \n {\bf Proof of Theorem 2.} We adapt an 
argument well known in the "local theory" of Banach spaces.
By Lewis' version of Fritz John's theorem (cf.
 [P5] p.~28) there is an isomorphism $u\colon \ E\to OH_n$
such that
 $\pi_{2,oh}(u) = \sqrt n$ and $\pi^*_{2,oh}(u^{-1}) = \sqrt n$. It is
 rather easy to check (cf.[P3]) directly from the
definition of the norm $\pi_{2,oh}$
 that for all $v\colon \ OH_n \to E$
 $$\pi^*_{2,oh}(v) = \inf\{\|B\|_{HS} \|A\|_{cb}\}$$
 where $B\colon\ OH_n \to OH_n$, $A\colon \ OH_n \to E$ and $v = AB$.
 
 \n Hence $u^{-1} = AB$ with $\|A\|_{cb}= 1$ and
$\|B\|_{HS}= \sqrt n$.
 Clearly $\|uA\|_{HS}\le \sqrt n$ by definition of $\pi_{2,oh}$, hence
 $$\|B\|_{HS} \ \|uA\|_{HS} \le n = tr(uu^{-1}) = tr(uA\cdot B)$$
 so by the equality case of the Cauchy Schwarz inequality we must have
 $$(uA) = B^*$$
 hence $B^{-1} = B^*$, so that $B$ is unitary. It follows that
 $$ {\|u^{-1}\|_{cb} \le \|A\|_{cb} \|B\|_{cb} 
 \le 1}$$
 since for $B\colon \ OH_n\to OH_n$ we clearly have
(cf.[P1-2]) $\|B\|_{cb}\le \|B\|$.
Conversely we have  $\sqrt n\, =
\pi_{2,oh}^*(u^{-1}) \le \sqrt n\,
 \|u^{-1}\|_{cb}$, which proves that
$\|u^{-1}\|_{cb}=1$.
 
\qed\medskip

 \n {\bf Proof of Corollary 3.} By Theorem 2 and by (5)
we have 
$$\pi^n_{2,oh }(I_E)^2\geq \pi_{2,oh}(I_E)^2/2,$$
from which the corollary follows.\qed

 \n {\bf Proof of Corollary 4.} By (6) we have
 $$\|u\|_{cb} \le \pi_{2,oh}(u)$$
 hence $\|u\|_{cb} \|u^{-1}\|_{cb} \le \sqrt n$.
 
 \n {\bf Proof of Corollary 5.} Let $u\colon \ E\to OH_n$
be as in
 Theorem~2. By Proposition~1 there is an extension $\tilde u\colon \ B(H)\to
 OH_n$ such that $\pi_{2,oh}(\tilde u) \le \sqrt n$. By
(6) we have $\|\tilde u\|_{cb} \le \sqrt n$, hence letting
$P =
 u^{-1}\tilde u$ we find $\|P\|_{cb} \le \sqrt n$.\qed\medskip
 
 Finally we have (by going through $OH_n$).
 
 \proclaim Corollary 6. Let $E,F$ be arbitrary
$n$-dimensional operator
 spaces. There is an isomorphism $u\colon \ E\to F$ such that $\|u\|_{cb}\
 \|u^{-1}\|_{cb}\le n$.
 
 Note that this is optimal (asymptotically) already in the category of
 Banach	spaces, as shown by the well known spaces constructed by
E.Gluskin
 [Gl1,2]. We refer the reader to [Pa2] for a
discussion of the problem considered in corollary 6 when
$E$ and $F$ are the same underlying Banach space. Even
when the  Banach space underlying $E$ and $F$ is the
$n$-dimensional Euclidean space, the asymptotic order of
growth of corollary 6 cannot be improved (see
Theorem 2.15 in [Pa2]).

\magnification\magstep1
\baselineskip = 18pt
\def\n{\noindent}
 
\overfullrule = 0pt

Finally we turn to some examples.\medskip
\item{1)} If $E_n = OH_n$, then clearly $\pi_{2,0h}(I_{E_n}) =
\pi_2(I_{E_n}) = \sqrt n$. More generally since we have a completely
contractive inclusion (cf [P1]) $OH_n \to R_n+C_n$, we have
$\pi_{2,oh}(I_{R_n+C_n}) = \sqrt n$.
\item{2)} If $E_n = R_n$ or $C_n$, we claim that
$$\pi_{2,oh}(I_{R_n}) =\pi_{2,oh}(I_{C_n}) = n^{1/4}.$$
Indeed, let $e_{1i}, \ldots, e_{1n}$ be the canonical basis of $R_n$. It is
easy to check that\break $\left\|\sum\limits^n_1 e_{1i} \otimes \bar
e_{1i}\right\|^{1/2} = n^{1/4}$. Since $\left(\sum
\|e_{1i}\|^2\right)^{1/2} = n^{1/2}$ we find $\pi_{2,oh}(I_{R_n}) \ge
n^{1/4}$. On the other hand, by the interpolation theorem
 for operator spaces (cf.[P1] Remark 2.11), for all
$u\colon \ OH\to R_n$ we have $\|u\|_{cb} = (tr|u|^4)^{1/4}$ where $|u| =
(u^*u)^{1/2}$ is the modulus of $u$ as an operator between Hilbert spaces.
Hence if $(T_i)$ denotes an othonormal basis of $OH$, we have since $R_n$
is $n$-dimensional
$$\eqalign{\left(\sum \|u(T_i)\|^2\right)^{1/2} &= (tr|u|^2)^{1/2} \le
n^{1/4}(tr|u|^4)^{1/4}\cr
&\le n^{1/4} \|u\|_{cb}}$$
hence $\pi_{2,oh}(I_{R_n}) \le n^{1/4}$.

This proves the above claim for $R_n$. The proof for $C_n$ is similar.
\item{3)} Let $u_1,\ldots, u_n$ be unitary operators in $B(H)$ such that
$u_i =u^*_i$, $u^2_i=I$ and
$$u_iu_j+ u_ju_i = 0 \quad \hbox{if}\quad i\ne j.$$
These are the canonical generators of a Clifford algebra. It is known that
such operators can be constructed inside the space $M_{2^n}$. Let $E_n =
\hbox{span}(u_1,\ldots, u_n)$. We claim that $\pi_{2,oh}(I_{E_n}) \le \sqrt
2$. Let $i\colon \ \ell^n_2\to M_{2^n}$ be the map defined by $i(x) = \sum
x_iu_i$.

Clearly we have
$$\forall\ x\in \ell^n_2 \qquad i(x)^* i(x) + i(x) i(x)^* = 2 \|x\|^2
I.\leqno  {(8)}$$
This implies
$$\forall\ x\in \ell^n_2 \qquad \|x\|^2 \le \|i(x)\|^2 \le
2\|x\|^2.\leqno  {(9)}$$ Moreover let us denote
$$\forall\ a \in M_{2^n}\qquad  \tau(a) = 2^{-n}tr(a).$$
Then the identity (8) yields
$$\tau(i(x)^* i(x)) = \|x\|^2$$
hence we have by (9)
$$\forall\ x\in \ell^n_2 \qquad {1\over 2}\|i(x)\|^2 \le
\|x\|^2 = \tau(i(x)^* i(x)) \le \|i(x)\|^2.\leqno
 {(10)}$$ Let us denote $\|a\|_2 = (tr\ a^*a)^{1/2}$
for all $a$ in $M_{2^n}$. Now consider a finite sequence
$(a_j)$ in $E_n$. Let $a_j = i(x_j)$ with $x_j \in
\ell^n_2$. We have $$\eqalignno{\left\|\sum a_j\otimes
\bar a_j\right\| &= \sup_{\|y\|_2\le 1} \left\|\sum
a_jya^*_j\right\|_2\cr &\ge 2^{-n/2} \left\|\sum
a_ja^*_j\right\|_2 \ge 2^{-n} tr\left(\sum a_ja^*_j\right)
= \sum \tau(a_ja^*_j)\cr \noalign{\hbox{hence by (10)}}
&\ge 2^{-1} \sum \|a_j\|^2.}$$
Hence we have $\pi_{2,oh}(I_{E_n}) \le 2^{1/2}$.

In [ER7], Effros and Ruan proved an analogue of the
Dvoretzky-Rogers theorem for operator spaces. We can deduce
a similar (and somewhat more precise) result from the
above  corollary 5. Indeed, let $E\subset B(H)$ be any
$n$-dimensional operator space. Let us denote by
$i_E\colon \ E\to B(H)$ the embedding. Then the ``operator
space nuclear norm'' of $i_E$, denoted by 
$\nu(i_\varepsilon)$, as introduced in [ER6] satisfies
$$\nu(i_E) \ge n^{1/2}.$$ Indeed, by corollary 5 there is a
projection $P\colon \ B(H)\to E$ with $\|P\|_{cb} \le
n^{1/2}$, hence $$\eqalign{n = tr(I_E) &= tr(Pi_E) \le
\|P\|_{cb} \nu(i_E)\cr &\le n^{1/2}  \nu(i_E).}$$ This
implies that the identity of an operator space $X$ is
``operator 1-summing'' in the sense of [ER7] iff $X$ is
finite dimensional.

\n {\bf Acknowledgement:} I am grateful to Vern Paulsen
for stimulating conversations.
 \vfill \eject

\centerline{\bf References}

 \item{[BP]} D. Blecher and V. Paulsen. Tensor products of operator spaces.
 J. Funct. Anal. 99 (1991) 262-292.

 \item{[B1]} D. Blecher. Tensor products of operator spaces II. (Preprint)
 1990. To appear in Canadian J. Math.
 
 \item{[B2]} D. Blecher. The standard dual of an operator
space.
 To appear in Pacific J. Math.

 \item{[ER1]} E. Effros and Z.J. Ruan. On matricially normed spaces. Pacific
 J. Math. 132 (1988) 243-264.
 
 \item{[ER2]} $\underline{\hskip1.5in}$. A new approach to operators spaces.
 (Preprint).

 \item{[ER3]} $\underline{\hskip1.5in}$. On the abstract characterization of
 operator spaces. (Preprint)
 
 \item{[ER4]} $\underline{\hskip1.5in}$. Self duality for the Haagerup
 tensor product and Hilbert space factorization.
 (Preprint 1990)
 
 \item{[ER5]} $\underline{\hskip1.5in}$. Recent development in operator
 spaces. (Preprint)
 
 \item{[ER6]} $\underline{\hskip1.5in}$. Mapping spaces
and liftings for operator spaces. (Preprint)

\item{[ER7]}  $\underline{\hskip1.5in}$. The
Grothendieck-Pietsch and Dvoretzky-Rogers Theorems for
operator spaces. (Preprint 1991)

\item{[Gl1]} E.Gluskin. The diameter of the Minkowski
compactum is roughly equal to $n$. Funct. Anal. Appl. 15
(1981) 72-73.

\item{[Gl2]} $\underline{\hskip1.5in}$. Probability in
the geometry of Banach spaces. Proceedings I.C.M.
Berkeley, U.S.A.,1986 Vol 2, 924-938 (Russian),
Translated in "Ten papers at the I.C.M. Berkeley 1986"
Translations of the A.M.S.(1990)

 \item{[Pa1]}$\underline{\hskip1.5in}$.  Completely
bounded maps and dilations. Pitman Research Notes 146.
Pitman Longman (Wiley) 1986.
 
\item{[Pa2]} V. Paulsen. Representation of Function
algebras, Abstract
 operator spaces and Banach space Geometry. J.
Funct. Anal. To appear.

 \item{[P1]} G. Pisier. The Operator Hilbert space and
complex interpolation. Preprint.

 \item{[P2]}$\underline{\hskip1.5in}$.  Espace de Hilbert
d'op\'erateurs et Interpolation complexe. C.R. Acad. Sci.
Paris To appear.
 
 \item{[P3]}$\underline{\hskip1.5in}$. Sur les op\'erateurs
factorisables par $OH$.  C.R. Acad. Sci.
Paris. To appear.

\item {[P4]} $\underline{\hskip1.5in}$. Factorization of linear
operators and the Geometry of Banach spaces.  CBMS
(Regional conferences of the A.M.S.)    60, (1986),
Reprinted with corrections 1987.

\item {[P5]} $\underline{\hskip1.5in}$. The volume of
convex bodies and Banach space geometry.
Cambridge Univ. Press, 1989.
 
\item{[Ru]} Z.J. Ruan. Subspaces of $C^*$-algebras. J. Funct. Anal. 76
 (1988) 217-230.
 
 \item{[TJ]} N.Tomczak-Jaegermann. Banach-Mazur distances
and finite dimensional operator ideals. Pitman,1988.

\vskip12pt\vskip12pt
Texas A. and M. University

College Station, TX 77843, U. S. A.

and

Universit\'e Paris 6

Equipe d'Analyse, Bo\^\i te 186,
 
75230 Paris Cedex 05, France

\end

 \end